\documentclass[review,hidelinks]{siamart171218}

\usepackage[utf8]{inputenc}

\usepackage{amsfonts}
\usepackage{dsfont}  %
\usepackage{mathrsfs} 
\usepackage{graphicx}
\usepackage{comment}
\usepackage{graphicx}      

\usepackage{amsmath,amssymb}
 \usepackage{latexsym}
 \usepackage{multirow}

 \usepackage{bbm}

\usepackage{array}
\usepackage{makecell}
\usepackage{multirow}
\usepackage{xcolor}
\usepackage[T1]{fontenc}
\usepackage{eulervm}

\usepackage{hyperref} 
\usepackage{tikz}

\makeatletter
\DeclareRobustCommand{\sqcdot}{\mathbin{\mathpalette\morphic@sqcdot\relax}}
\newcommand{\morphic@sqcdot}[2]{%
  \sbox\z@{$\m@th#1\centerdot$}%
  \ht\z@=.33333\ht\z@
  \vcenter{\box\z@}%
}
\makeatother

\definecolor{purple}{RGB}{69,22,170}

\definecolor{verde}{RGB}{83,134,31}

\newcommand{\trp}{\intercal}

 \DeclareMathOperator{\Diag}{Diag}
 \DeclareMathOperator{\diag}{diag}
 \DeclareMathOperator{\tr}{tr}
 \DeclareMathOperator{\sign}{sign}
 \DeclareMathOperator{\cov}{Cov}

\ifpdf
  \DeclareGraphicsExtensions{.eps,.pdf,.png,.jpg}
\else
  \DeclareGraphicsExtensions{.eps}
\fi

\newsiamremark{assumptions}{Assumptions}

\newsiamremark{remark}{Remark}

\ifpdf
\hypersetup{
  pdftitle={The $\mathcal{H}_2$-optimal Control Problem of CSVIU Systems: Discounted, Counter-discounted and Long-Run Solutions~Part I: The Norm},
  pdfauthor={Jo\~ao B. R. do Val and Daniel S. Campos}
}
\fi

\title{The $\boldsymbol{\mathcal{H}_2}$-optimal Control Problem of CSVIU Systems: Discounted, Counter-discounted and Long-Run Solutions\\Part I: The Norm%
\thanks{Submitted to the SIAM J Control  Optim editors on \today
\funding{Research supported in part by the National Council for Scientific and Technological Development (CNPq), grant n. 303352/2018-3, by FAPESP under grant n. 2016/08645-9, and by the Coordenação de Aperfeiçoamento de Pessoal de Nivel Superior - Brasil (CAPES) - Finance Code 001.}}}

\author{ Jo\~ao B. R. do Val\thanks{
UNICAMP, School of Electrical and Computer Engineering, Av. Albert Einstein 400, Campinas, SP, Brazil (\email{jbosco@fee.unicamp.br}, \email{d211498@dac.unicamp.br})}
\and
Daniel S. Campos\footnotemark[2]
}

\headers{The $H_2$-optimal Control of CSVIU Systems -- Part I}{Jo\~ao B. R. do Val and Daniel S. Campos}

\begin{document}
\bibliographystyle{siamplain}

\maketitle

\begin{abstract}                
The paper deals with the $H_2$-norm and associated energy or power measurements for a class of processes known as CSVIU (Control and State Variation Increase Uncertainty). These are system models for which a stochastic process conveys the underlying uncertainties, and are able to give rise to cautious controls. The paper delves into the non-controlled version and fundamental system and norms notions associated with stochastic stability and mean-square convergence. 
One pillar of the study is the connection between the finiteness of one of these norms or a limited energy measurement growth with the corresponding stochastic stability notions. A detectability concept ties these notions, and the analysis of linear-positive operators plays a fundamental role. 
The introduction of various $H_2$-norms and energy measurement performance criteria allows one to span the focus from transient to long-run behavior. As the discount parameter turns into a counter-discount, the criteria enforce stricter requirements on the second-moment steady state errors and on the exponential convergence rate. A tidy connection among this $H_2$-performance measures cast employs a unifying vanishing discount reasoning.

\end{abstract}



\section{Introduction}

An exemplary mathematical model of a process of interest forms the ideal basis for forecasting future behavior, evaluating critical phenomena or measurements, or designing a control system.  Nevertheless, models are mostly nothing more than feeble imitations of reality, even though a good model can carry many of the essential features of the real world, \cite{bib:1,MKac}. As pointed out by many authors, this idea translates, in part, that to have more representative models for existing systems, one ought to characterize the uncertainties adequately,  c.f.  \cite{FrancisKhar,DraganMorozanStoica,bib:2}.  The need to deal with poor models is one of the main drives in the control field during the last three decades, which has brought a flourishing of new ideas that permeates many of the present date research.

Along these lines, the authors claim that the CSVIU model is a stochastic-based alternative to the uncertainty representation inherent to poor modeling of dynamic systems, cf.  \cite{bib:4,Fernandes2020,bib:3}. The CSVIU concept avoids the usual worst-case analysis of robust control by introducing a particular stochastic perturbation. A feature of the CSVIU method is a perturbing stochastic process that produces more significant drifts as the system state and control deviate from a better-known operation point, providing an interesting ground for accounting modeling errors.  The approach builds on the idea that when models are frails, state and control variations increase uncertainty, and cautious controls should prevail. The idea is to inbuild mathematically the notion that conservativeness should take place in the face of uncertainties. Solving the underlying optimal stochastic control problem, a region on the state space arises, in which an inaction control is optimal.  The zero-control (or zero-variation control) is the optimal feedback control therein, a behavior not encountered in the robust worst-case analysis, as far as the authors are aware.

The present paper focuses on the $H_2$-norm, or energy and power inspired measurements in associated settings. The companion paper \cite{CSVIU:control} presents the optimal $H_2$-norm and the optimal overtaking control, relying on the calculus developed here. Let us consider the following discrete-time process defined in a complete filtered probability space, $(\Omega, \mathcal{F}, P, \{ \mathcal{F}_k \}_{k \geq 0})$, devised by the CSVIU approach,
\begin{equation}
\Theta := \left\{
\begin{aligned}
    &x(k+1) = A x(k) + ( \sigma_x + \overline{\sigma}_x \diag(|x(k)|)) \varepsilon (k) + \sigma \omega (k), \\
&y(k) = Cx(k),  \quad x(0) = x
\end{aligned}\right.
\label{CSVIUsystem_no_control}
 \end{equation}
 where $x\in\mathds{R}^n$, $A, \sigma_x, \overline{\sigma}_x \in \mathds{R}^{n \times n}$, $\sigma\in\mathds{R}^{n\times r}$ and $C\in \mathds{R}^{p \times n}$. The noise $\{\omega (k)\}_{k\ge0}$ is a $r$-dimensional persistent disturbance noise sequence, the ``nature noise'', and the extra noise $\{ \varepsilon (k)\}_{k \geq 0}$ is the  $n$-dimensional ``intrinsic model noise'' sequence. Both are i.i.d. sequences with zero mean and their joint covariance forms an identity matrix. The filtration $\{ \mathcal{F}_k \}_{k \geq 0}$ is composed by the sub-$\sigma$-algebras $\mathcal{F}_k \subset \mathcal{F}$ generated by the random variables $\omega(0),\varepsilon(0), \ldots, \omega(k),\varepsilon(k)$. Given a $n$-valued vector $x, \diag(|x|)$ is the $n$-dimensional diagonal matrix formed by setting $|x| = [\begin{smallmatrix} |x_1|& |x_2|&\cdots & |x_n|\end{smallmatrix}]^\trp$ as its diagonal, where $| \cdot |$ is the absolute value. The processes $\{x (k)\}_{k\ge0}$ is the $n$-dimensional state and $\{y (k)\}_{k\ge0}$ is a $m$-dimensional output of interest with $m\le n$.

The model $\Theta$ is an attempt to describe the unknown real system, 
\begin{equation}
\begin{aligned}
{z(k+1)}& = f(z(k)) + \sigma \omega(k)
\end{aligned}
\label{autonomus_system}
\end{equation}
for some  $f:\mathds{R}^n\to\mathds{R}^n$, near a point $\bar{z}$. In this account, the dynamic matrix $A$ in  \eqref{CSVIUsystem_no_control} represents the Jacobian matrix of $f$ at $\bar{z}$, assuming that the derivatives exist. Since $f$ is a not well-known system function, matrix $A$ is bound to be a poor representation of the Jacobian. For such a linear representation one can write that 
\begin{equation}
x(k+1) = Ax(k) + f(\bar{z})-\bar{z} 
+  \Bigl(\frac{\partial f}{\partial z}\Big|_{\bar{z}}-A\Bigr)x(k) + o(|z(k)-\bar{z}|^2)+ \sigma \omega(k)
\label{autonomus_system.var}
\end{equation}
 holds, where we set $x(k) = z(k) -\bar{z}$. 
 
 If $\bar{z}$ is a precisely known equilibrium point, the first difference on the rhs is null, and also null is the second difference, if the Jacobian of $f$ at $\bar{z}$ is precisely known. In an attempt to represent these residuals, the CSVIU  model $\Theta$ includes the extra-terms associated with the noise sequence $\{\varepsilon(k)\}_{k\ge0}$. 
The additional noise terms represent the error of such an educated guess;  ${\sigma}_x\varepsilon(k)$ mainly for the mismatch $\delta f = f(\bar{z})-\bar{z}$, whereas  $\bar{\sigma}_x\diag(|x_t|)\varepsilon(k)$ stands for the lack of trust on the linearized model, as each of the state vector component displaces from the chosen point $\bar{z}$ (or $x=0$). 
The errors due to higher orders terms are expressed by the residue vector function $o(\cdot)$, depending on the square of the componentwise distances $|z_i(k)-\bar{z}|, i=1,\dots,n$ for which the term in \eqref{CSVIUsystem_no_control} involving $\diag(|x|)$ tries to convey in terms of larger variance as such displacements increase.

The representation of linear systems in continuous-time driven by a nonlinear diffusion term of Brownian Motion appears in the stochastic literature, known as perturbed linear systems, or as semilinear SDEs, e.g., \cite[cap\;4]{Mao.1991} or  \cite[cap\;6.7]{Khasminskii.2012}.  It is motivated, among other reasons, to study a nonlinear stochastic system subject to poor modeling. The linear part of the model comes as the best guess of the Jacobian matrix evaluated at some point of interest or the point of maximum knowledge of the original system. The CSVIU model, such as $\Theta$, emerges as a specific discrete-time counterpart for those system models.  

This paper studies the discounted and non-discounted $H_2$-norm type of performance problems of non-controlled CSVIU systems. It tailors appropriated notions of stochastic stability to the norms or energy measurements proposed here. It makes these notions meaningful by connecting the finiteness of these norms with a corresponding idea of stability. In some situations, the $H_2$-norm may be unbounded due to the persistent noise. Nevertheless, this scenario can be dealt with in a problem for which the energy increases at a limited rate, giving origin to the notion of an overtaking measurement induced by such a norm-like measurement.

For a $\{ \mathcal{F}_k\}_{k\ge0}$-adapted $n$-dimensional process, $k \to w(k)$, $k\ge0$ defined on a filtered probability space $( \Omega, \mathcal{F}, P, \{  \mathcal{F} \}_{k \geq 0})$, let us consider the following $\ell_2 (\Omega, \mathcal{F}, P)$ mean energy measurements.  Set for $\kappa\ge0$ and some $\alpha>0$,
\[
\mathcal{E}^{\kappa,\alpha}_2(w(\cdotp)) := E_{x} \Bigl[  \sum_{k=0}^{\kappa} \alpha^k \| w(k) \|^2 \Bigr]
\]
in which $E_{x} \big[ \cdot \big]$ is the short for the expectation $E \big[ \cdot | \mathcal{F}_0 \big]$. The $H_2$-$\alpha$-norm of $w(\cdot)$ is defined as
\begin{equation}
    \mathcal{E}^\alpha_2(w(\cdotp)):=\lim_{\kappa\to\infty}\mathcal{E}^{\kappa,\alpha}_2(w(\cdotp)), \quad w(0)=0, 
    \label{energymeasure1}
\end{equation}    
whenever finite. The energy measurement is in Abel's mean form, and if $\alpha <1$, it is a discounted measurement, which is adequate to deal with the persistent noise excitation in the $\Theta$ system. 
 
 When $\alpha\ge1$, we called it a counter-discounted measurement, the limit in  \eqref{energymeasure1} is possibly unbounded, and the norm is not defined.  But we can still measure the energy of signals and compare them against each other by means of an overtaking criterion. Namely, two $\{ \mathcal{F}_k\}_{k\ge0}$-adapted $n$-dimensional processes, $k \to w(k)$ and  $k \to z(k)$ are comparable by overtaking if for each $\kappa> \kappa_0$,
 \begin{equation}\label{eq:overtaking.1}
     \mathcal{E}^{\kappa,\alpha}_2(w(\cdotp))\le \mathcal{E}^{\kappa,\alpha}_2(z(\cdotp))+\epsilon
 \end{equation}
 
 A third measurement of interest is based on the C\`esaro's mean form,
\begin{equation}
    \mathcal{P}_2(w(\cdotp)) := \lim_{\kappa \to \infty} \frac{1}{\kappa} E_{x} \Bigl[  \sum_{k=0}^{\kappa-1} \| w(k) \|^2 \Bigr]
    \label{energymeasure2}
\end{equation}
provided that the limit exists.  Possible $H_2$-norms of interest for $\Theta$ are  $\mathcal{E}^\alpha_2(y(\cdotp))$ with $x(0)=0$ and $\mathcal{P}_2(y(\cdotp))$. Throughout they are referred as $\mathcal{E}^\alpha_2$ and $\mathcal{P}_2$ for short.

The average power of the process $w(\cdotp)$ in \eqref{energymeasure2} is a notion connected with recurrence and finiteness of mean recurrence time to a compact set, cf.  \cite{MeynTweedie,Walters}. It is linked to the existence of a stationary distribution and ergodic behavior, as studied, for example, in \cite{bib:11}, concerning a different class of processes.
 
 Specific norms are an essential form of rendering performance and robustness measurements for a control system in the study of signals and systems. Usually, the appropriate norm depends on the situation at hand,  associated somehow with the energy or error of the system output. The conception of norms for deterministic systems has reached maturity, appearing on various studies in the control literature spanning for more than three decades.  Note the interesting fact that apart from a linear system driven by additive noise, deterministic $H_2$ and $H_\infty$ norms are never be equivalent to their stochastic counterparts. 
 
 A deterministic approach to the norm of stochastic systems does not make sense. Besides mere inadequacies, the usual notions are not applicable due to the noise persistence, driving system variables into permanent fluctuations.
 Regarding the definition and computation of norms for stochastic systems, there exist notions such as studied in \cite{bib:2,bib:11,bib:17,bib:18,bib:19}, for certain classes of system. However, there lies much room for stochastic suitableness.

The $\alpha$-discounted measure in Abel's mean form of the expected energy in \eqref{energymeasure1} stands as a possible stochastic norm. 
The energy measurement in \eqref{energymeasure1} when $\alpha<1$ attaches a discount factor to the future errors, an approach that focuses on system energy expenditure on the near horizon, thus, on the transient behavior of the system. 
The power measure in C\`{e}saro's mean form in \eqref{energymeasure2} represents the average mean norm, and it connects with the study of the long-run average error appearing in many contexts. 
The $\alpha$-counter-discounted measurement, as a rule, will be unbounded, but the expected energy still can be compared by the overtaking criterion.

A further step is to unify the treatment of these classes of indices based on norm criteria. The Abel and the C\`esaro mean can be related by limiting reasoning known as \emph{vanishing discount}, and we explore this connection here. Tauberian theorems play a vital role to deal with the long-run average $H_2$-norm, cf. \cite{bib:6}.  As a result, most of the relations concerning the discounted/counter-discounted problems also hold for the average problem with $\alpha=1$.

   \paragraph{Stochastic Stability}
 
 Essential to dynamical systems is the stability notion, here understood in proper stochastic senses. Consider the  following notions.
 
\begin{definition}
    System $\Theta$ is \emph{$\alpha$-stochastic stable} if,
    \begin{enumerate}
        \item[i)] $0\le\alpha<1$, and the measurement $\mathcal{E}^\alpha_2(x(\cdotp)) < \infty, \; \forall x(0) = x_0 \in \mathds{R}^n$, 
        
        \item[ii)] $\alpha\ge1$, and there exist $c_0,c_1$  and $\xi_\kappa\in\mathds{R}^n$ with $\|\xi_\kappa\|\le c_2$ for each $\kappa>0$  
        such that,  
        \[
        \mathcal{E}^{\kappa,\alpha}_2(x(\cdotp)) \le c_0\|x_0-\xi_\kappa\|^2+ \kappa c_1\alpha^\kappa,  \;\forall x(0) = x_0 \in \mathds{R}^n.
        \]

    \end{enumerate}
    \label{stability_definition1}
\end{definition}

When $\alpha =1$, stability is especially highlighted in the C\`{e}saro mean form.
\begin{definition}
    System $\Theta$ is \emph{stochastic stable} if the measurement $\mathcal{P}_2(x(\cdotp))\le \bar{c} < \infty$ for any $x(0) = x_0 \in \mathds{R}^n$.
    \label{stability_definition2}
\end{definition}

Commonly, some specified system energy, such as that carried by the output $y(\cdot)$, expresses the system's performance. At the same time, we seek stability in one of the senses above.  However, the relation between such a measurement of interest and the system stability is, in general, not known. More specifically, how finiteness of an energy measurement can assure stability?
When the norms $\mathcal{E}^\alpha_2$ or $\mathcal{E}_2$ are finite and well defined, do these imply stochastic stability for system $\Theta$ in the corresponding senses of Definitions \ref{stability_definition1} or \ref{stability_definition2}? 
A positive answer comes from a link between the finiteness of such measurements with the respective notions of $\alpha$-stability and stochastic stability.  It turns out that proper concepts of detectability tailored for the class of CSVIU models set up these bridges in an undeniable way. For that, the analysis of linear-positive operators plays a fundamental role associated with perturbed Lyapunov equations.

After this critical step, stability analysis comes in hand, and one can deal with the discounted norm quite straightforwardly.  However, this option leads to a stability requirement that may often be too loose since $\alpha$-stability only implies that $\alpha^{k/2}\|x_k\|^2 \to 0$ in the mean-square sense, cf. \cite{bib:4,bib:3}.

On the other hand, the counter-discounted ($\alpha\ge1$) and the long-run average approaches are more defiant and entirely new for CSVIU models. Adopting the measurements in \eqref{eq:overtaking.1} or \eqref{energymeasure2}, the analysis encompasses the system performance in the far future (or on the steady state), taking a better view into the asymptotic behavior of state and output processes.  When $\alpha>1$, the overtaking criterion even more severely scrutinizes the steady state behavior and brings forth stronger stability conditions in the counter-discounted formulation. Complementary to the discounted and average $H_2$-norms, one can assess the convergence to a value in the mean-square sense and the convergence rate in exponential form.

This energy measurement span allows one to cover specific interest problems, from a focus on the transient behavior to the steady state, through different and more stringent requirements on the system's stochastic stability, as one increases the parameter $\alpha$.  The long-run power average and the overtaking energy criteria proceed neatly from the results developed here.  With these elements, we depart from the discounted LQ-control problem and weaken the requirements for the stability and stability notions studied in \cite{bib:3} to a broader class of $H_2$-norm like problems and ensuing stability notions.

 The paper reads as follows. Section \ref{sec:energyandmeasurement} develops in Lemma \ref{lemma:time.forward} the calculus on expected evolutions  and connects the notions of stability with the study of linear-positive operators in Proposition \ref{prop:freiling.stab}. The requirements for the association between finiteness or limited growth of energy measures \eqref{energymeasure1} and \eqref{energymeasure2} and stability of $\Theta$ is weakened in Theorem \ref{thm:detectability} with the detectability notion.
 
 Section \ref{subsec:long-run.ana} completes the relation between energy and power measurements, and $H_2$-norms conceived here. In Lemma \ref{lemm:convergence.alpha.ge1} notions of the state asymptotic behavior in the mean-square sense when parameter $\alpha\ge1$ appear, see also Remark \ref{rem:convergence.properties}. Theorem \ref{thm:analysis:longrunnorm2} develops the link between energy and power measurements studied. Finally, $H_2$-norms, energy measurements and associated convergence notions are dealt with in Corollaries \ref{coro:discounted.norm} and \ref{coro:average.norm}. Section \ref{sec:conclusion} brings some concluding remarks.

\section{Energy Measurements and  Stability}
\label{sec:energyandmeasurement}

For a matrix $U\in\mathds{R}^{n\times n}$, $\Diag(U) \in \mathds{R}^{n \times n}$ is the diagonal matrix formed by the main diagonal of $U$, $U_d\in\mathds{R}^n$ is the vector in the main diagonal of $U$ and $\tr(U)$ denotes the trace operator.  $\mathds{S}^{n+}$ stands for the real vector space of symmetric matrices of size $n$ that are positive semidefinite. For  $U\in\mathds{S}^{n+}$, $U \succ 0$ ($U\succeq0$) indicates that $U$ is a positive (semi-) definite matrix, also, for $U,V\in\mathds{S}^{n+}$, $V\succeq U \Leftrightarrow V-U\succeq 0$. $\|U\|$ indicates any matrix norm and for a square matrix $\lambda^+(U) (\lambda^-(U))$ denotes the largest (smallest) eigenvalue of $U$ and $r_\sigma(U)$ its spectral radius.
 
For a vector $u \in \mathds{R}^n$, $|u|$ indicates the vector $[\begin{smallmatrix}
|u_1|&|u_2|&\cdots&|u_n|
\end{smallmatrix}]^\trp$ and $\diag(u)$ stands for the diagonal matrix made up by vector $u$.
For two vectors $u,v \in \mathds{R}^n$, $\langle u , v \rangle$ denotes the inner product,  $u \sqcdot v$ denotes the Hadamard product, and the square (semi-)norns $\|x\|^2_U$ stands for  $\langle x, U x\rangle$, in which $U \in \mathds{S}^{n+}$.

In connection with the data of system $\Theta$, let us define the linear operators $\mathcal{Z}: \mathds{S}^{n+} \to\mathds{S}^{n+}$, $\mathcal{W}: \mathds{S}^{n+} \to \mathds{S}^{n}$ and $ \varpi : \mathds{S}^{n+} \to \mathds{R}$, given by:
\begin{subequations}\label{eq:state-operators}
\begin{align}\allowdisplaybreaks
    \mathcal{Z}(U) &= \Diag ( \overline{\sigma}_{x}^{\trp} U \overline{\sigma}_x ),
    \label{linear_operators_definition}
\\
\mathcal{W}(U)  &= \Diag(\overline{\sigma}_x^\trp U \sigma_x + \sigma_x^\trp U \overline{\sigma}_x ),
     \label{linear_operators_definition.2}
\\
    \varpi(U)  &= \tr \{ U(\sigma \sigma^\trp + \sigma_x \sigma_x^\trp) \},
\intertext{together with $\mathcal{L}^\alpha:\mathds{S}^{n+} \to \mathds{S}^{n+}$, for some $\alpha \geq 0$,}
    \mathcal{L}^\alpha(U) &= \alpha(A^\trp U A + \mathcal{Z}(U)), 
    \label{special_linear_operators_definition}
\end{align}
\end{subequations}
and sometimes we denote $\mathcal{L}^1(U) = \mathcal{L}(U)$.
With the exception of $\mathcal{W}$, these are linear-positive operators, i.e., $U \succeq 0$ implies $\Pi(U) \succeq 0$, where $\Pi$ stands for any of the operators in \eqref{eq:state-operators}. 
%
In addition, they are monotone operators, since if $U \succeq V$ for $U,V \in \mathds{S}^{n+}$ then $\Pi(U)\succeq\Pi(V)$. 
When convenient, we refer by $\mathcal{W}_d(U)$ to the $\mathds{R}^n$-vector in the diagonal matrix  $\mathcal{W}(U)$.

Let us consider the set $\{ -1, 0, +1\}$ and create the collection $\mathscr{S}= \{ -1, 0, +1\}^n$ of distinct vectors $s_i \in\mathscr{S}, i= 1, \ldots , 3^n $ under some arbitrary order.
We introduce the signal vector function $\mathcal{S}: \mathds{R}^n \to \mathscr{S}$ such that for each $x\in\mathds{R}^n$,
\begin{equation}
   \mathcal{S}(x) = \begin{bmatrix}\sign(x_1) & \cdots & \sign(x_n)\end{bmatrix}^\trp
   \label{eq:signalfunction}
\end{equation}
with the convention $\sign(0)=0$.

For subsequent developments, let us denote for some $\alpha>0$, $Q\succeq0$ and $\kappa>0$,
\begin{equation}
    \mathcal{E}^{\alpha, \kappa}_{2, Q}(x(\cdotp)):=  E_x \Bigl[ \sum_{k=0}^\kappa \alpha^k \| x(k) \|_Q^2  \Bigr],\quad x(0)=x
    \label{discounted_finite_time_horizon}
\end{equation}
When $\kappa\to\infty$ the functional, possibily infinite, is called the \emph{$Q$-mean energy} of system $\Theta$. 
For further use, in the next lemma we consider a version of system $\Theta$ with exogenous input,
\begin{equation}
\Theta_{\textrm{ex}} := \left\{
\begin{aligned}
    &x(k+1) = A x(k) + B \ell(k)
 + ( \sigma_x + \overline{\sigma}_x \diag(|x(k)|)) \varepsilon (k) + \sigma \omega (k), \\
&y(k) = Cx(k) + D \ell(k),  \quad x(0) = x
\end{aligned}\right.
\label{CSVIUsystem_exogenous}
 \end{equation}
where $B\in \mathds{R}^{n \times m}$ and $D\in \mathds{R}^{p \times m}$. The $m$-dimensional process $k\to\ell_k$ is $\{\mathcal{F}_k\}_{k\ge0}$-adapted and it is such that $P(\ell_{n+k}\in A|\mathcal{F}_k) = P(\ell_{n+k}\in A|x_k),\forall A\in \mathcal{B}(\mathds{R}^m), n\ge0$.

Sometimes we opt for the more compact notation for system $\Theta$ or $\Theta_{\textrm{ex}}$,
\begin{equation}\label{eqs:compact.theta.1}
x_{k+1} = A x_k + B\ell_k + \sigma(x_k)\omega_{0,k}
 \end{equation}
 such that 
 $   \sigma(x)  = \begin{bmatrix}
     \sigma & \sigma_x + \bar{\sigma}_x \diag(|x|)
    \end{bmatrix}
$%
 and $\omega_{0,k} = [ \omega_k ~~ \varepsilon_k  ]^\trp$ (recall that $\cov(\omega_{0,k})=I_{r+n}$). The associated processes and time stage sequences are similarly expressed using the stage $k$ as a subindex.

The next lemma will be useful for framing the picture that so far we sketched.

 \begin{lemma}\label{lemma:time.forward}
Consider the sequences $P_k\in\mathds{S}^{n+}, v_k\in \mathds{R}^n$ and $g_k\in\mathds{R}$, $k=0,1,\ldots,\kappa$ that satisfy the following difference equations,
\begin{subequations} \label{algebraicequationsH2normcoef}
\begin{gather}
 \mathcal{L}^\alpha(P_{k+1}) + Q = P_k,
\label{recursivelyapunov}
\\
A^\trp  (v_{k+1}  + 2 P_{k+1} B\ell_k) + \mathcal{W}(P_{k+1})\mathcal{S}(x_k)  = \alpha^{-1} v_{k},
\label{recursivevector.2}
\\
g_{k+1} +  \varpi(P_{k+1})+ \| B \ell_k \|^2_{P_{k+1}} + \langle v_{k+1},   B \ell_k \rangle = \alpha^{-1}  g_k,
\label{recursiveconstant}
\end{gather}
\end{subequations}
with $P_{\kappa} = \Phi\in\mathds{S}^{n+}, v_{\kappa} = \theta\in\mathds{R}^n$ and $g_{\kappa} = \gamma\ge0$. Then,
\begin{equation}
     \mathcal{E}^{\alpha, \kappa-1}_{2,Q}(x(\cdotp)) = \| x \|^2_{P_{0}} +  \langle E_{x}[v_0], x  \rangle + g_{0}   \\
    - \alpha^{\kappa} E_{x} \bigl[ \| x({\kappa}) \|^2_{\Phi} + \langle \theta, |x({\kappa})|  \rangle + \gamma \bigr].
\label{rewriting_equation.2}
\end{equation}
whenever $x(0) = x$.
\end{lemma} 
 
\begin{proof}
Note first that the operator involved in \eqref{recursivelyapunov} is linear-positive,  hence, $P_k\in\mathds{S}^{n+}\;\forall k$. 
Set the function $\phi:\mathbb{N}\times \mathds{R}^n\to \mathds{R}$ made up from the sequences $P_k, r_k$ and $g_k$ as, 
\begin{equation*}
    \phi(k, x) := \alpha^k (x^\trp P_k x + \langle r_k , |x| \rangle + g_k),\quad x\in\mathds{R}^n.
\end{equation*}
 
Having in mind that for any $r,x \in \mathds{R}^n$,$\langle r, |x| \rangle = \langle \mathcal{S}(x), r \sqcdot x \rangle= \langle\mathcal{S}(x)\sqcdot r,  x \rangle$ and $U \in \mathds{S}^{n+}$,  $\tr \{ U  \diag(|x|)\}= \langle \mathcal{S}(x), \Diag(U) x \rangle$ hold true, and accounting for  the dynamics of system $\Theta_{\textrm{ex}}$, one can evaluate a variation of $\phi$ along a path $k\to x_k$ (with the compact notation in \eqref{eqs:compact.theta.1}),
\begin{multline} \label{variationofthenorm}\allowdisplaybreaks
          \alpha^{-k} \big( \phi(k+1, x_{k+1}) - \phi(k, x_k) \big)=
          \\ 
          \alpha \| x_{k+1}\|^2_{P_{k+1}} +  \alpha \langle s_{k+1}, r_{k+1} \sqcdot x_{k+1} \rangle  + \alpha g_{k+1} -
          (\| x_{k}\|^2_{P_{k}} + \langle s_k, r_{k} \sqcdot x_{k} \rangle + g_k )=
\\
\alpha \| A x_k + B \ell_k \|^2_{P_{k+1}}   +  2 \alpha (A x_k +B \ell_k)^\trp P_{k+1} \sigma(x_k) \omega_{0_k}
           + \alpha \| \sigma(x_k)\omega_{0_{k}} \|^2_{ P_{k+1} }  - \| x_k \|^2_{P_k } +
          \\
      \alpha \langle s_{k+1}, r_{k+1} \sqcdot \big( A x_k + B \ell_k +  \sigma(x_k)\omega_{0_k} \big) \rangle - \langle s_k, r_{k} \sqcdot x_{k}  \rangle +  
          \alpha g_{k+1} - g_{k} 
\end{multline}
where we set $s_{k} = \mathcal{S}(x_{k})$ and $s_{k+1} = \mathcal{S}(x_{k+1})$. Note that,
\begin{multline*}
     E  [\| \sigma(x_k)\omega_{0_k} \|^2_U|x_k =x]  = \tr \{ \Diag(U \sigma(x) \sigma(x)^\trp) \}=
     \\
     x^\trp \mathcal{Z}(U) x + \tr \{\mathcal{W} (U) \diag(|x|)\} + \varpi(U)=
     \| x \|^2_{\mathcal{Z}(U)} + \langle \mathcal{S}(x), \mathcal{W}(U)x \rangle + \varpi(U)
\end{multline*}

Returning to \eqref{variationofthenorm},
\begin{multline}   \label{identity_with_mk}
      \alpha^{-k} \big( \phi(k+1, x_{k+1}) - \phi(k, x_k) \big)=
     \alpha \big( \| A x_k \|^2_{P_{k+1}} + \| x_k \|^2_{\mathcal{Z}(P_{k+1})} \big)- \| x_k \|^2_{P_k } +
     \\
        \alpha \langle s_{k+1}, r_{k+1} \sqcdot( A x_k + B \ell_k)\rangle + \langle s_k, \left(\alpha \mathcal{W}_{d}(P_{k+1}) - r_{k}\right) \sqcdot x_{k} \rangle + 
    \\
     \alpha \left( g_{k+1} + \| B \ell_k \|^2_{P_{k+1}} + 2\langle  B \ell_k, P_{k+1}Ax_k  \rangle +  \varpi(P_{k+1}) \right) - g_{k} + m_k=
\end{multline}
\vspace{-5ex} 
\begin{multline*}
     \langle x_k, \big(\alpha A^\trp P_{k+1} A+ \alpha \mathcal{Z}(P_{k+1}) -  P_k \big) x_k\rangle +
     \\
        \langle \alpha A^\trp ( r_{k+1} \sqcdot s_{k+1} + 2 P_{k+1}B \ell_k)  + \alpha \mathcal{W}(P_{k+1})s_k - s_k\sqcdot r_{k}   ,x_{k} \rangle +
     \\
     \alpha \bigl( g_{k+1} + \| B \ell_k \|^2_{P_{k+1}}  + \langle r_{k+1} \sqcdot s_{k+1},   B \ell_k \rangle     + 
    \varpi(P_{k+1}) \bigr) - g_{k} + m_k
\end{multline*}
where, the random process $k\to m_k$ is
\begin{equation*}
         m_k  :=  2 \alpha (A x_k+B\ell_k)^\trp  P_{k+1} \sigma(x_k)\omega_{0_k} +
         \\
           \alpha \langle s_{k+1}, r_{k+1} \sqcdot \sigma(x_k)\omega_{0_k}  \rangle, 
\end{equation*}
 a zero $\{\mathcal{F}_k\}$-martingale that comprises each of the remaining terms of \eqref{identity_with_mk} such that $E[m_k|x_k=x]=0$. By setting
$v_k:=s_k\sqcdot r_{k}, k=0,\ldots,\kappa-1$ we get that if 
\begin{gather*}
   \alpha \left( A^\trp P_{k+1} A+ \mathcal{Z}(P_{k+1})\right) + Q -  P_k =0,
   \\
    \alpha \left(A^\trp (v_{k+1}  + 2 P_{k+1}B \ell_k) +  \mathcal{W}(P_{k+1})s_k\right) - v_{k} =0,
    \\
         \alpha \left( g_{k+1} + \| B \ell_k \|^2_{P_{k+1}} + \langle v_{k+1},   B \ell_k \rangle     +  \varpi(P_{k+1}) \right) - g_{k} =0,
 \end{gather*}
one has that
\begin{equation}
    E[\phi(k+1, x_{k+1}) - \phi(k, x_k) | x_k] = - \alpha^{k} E[\| x(k) \|_Q^2|x_k]
\end{equation}
and since the process $\Theta_{\mathrm{ex}}$ is Markovian,
\begin{equation}
    E[\phi(\kappa, x_{\kappa}) - \phi(0, x_0)] = - E[\sum_{k=0}^{\kappa-1} \alpha^k \| x(k) \|_Q^2]
    \label{rewriting_equation}
\end{equation}
with $\phi(\kappa,x_\kappa) = \alpha^\kappa (x_\kappa^\trp \Phi x_\kappa + \langle \theta , |x_\kappa| \rangle + \gamma)$. 
Finally, one can rewrite \eqref{rewriting_equation} as \eqref{rewriting_equation.2} with the conditions in the lemma.
\end{proof}

\subsection{Stochastic Stability}

Stochastic Stability relates to finiteness of the state energy measurements in Definitions \ref{stability_definition1} (i) and \ref{stability_definition2}, or to the limited geometric growth in Definition \ref{stability_definition1} (ii). We approach here these stability notions, providing conditions for them to hold. They involve a perturbed type of Lyapunov equation fundamentally, but that, in general, does not suffice, as we will see.

In the next steps, the matrix equation
\begin{equation}\label{lyapunovkind.1}
     (I-\mathcal{L}^\alpha)(U) = Q 
\end{equation}
for $\alpha\ge 0$ and $Q\succeq 0$ plays a vital role. These are studied in terms of positive operators in ordered Banach spaces, connected to the idea of linear perturbed Lyapunov equations that arise in many stochastic problems, e.g., \cite{freiling,Hasanov, Zhao-Yanetall}. We bring the essentials of this analysis, adapted to the setting here. Recall that we deal with the linear-positive operators in \eqref{eq:state-operators}, with the exception of $\mathcal{W}$.

\begin{proposition}\label{prop:freiling.stab}
 $\alpha$-stability of $\Theta$ for some $0\le\alpha<1$ is equivalent to require that
\begin{enumerate}
    \item[i)] $\mathcal{L}^\alpha$ is an inverse-positive operator,
    \item[ii)] $\mathcal{L}^\alpha$ is $d$-stable,
    \item[iii)] There exists $U\succ0$ such that $(I-\mathcal{L}^\alpha)(U)\succ0$,
    \item[iv)] $\alpha^{1/2}A$ is $d$-stable relative to $\alpha \mathcal{Z}$,
    \item[v)] All eigenvalues of $\sqrt{\alpha}A$ lie in the open unit disk and $r_\sigma((I-\alpha\mathds{A})^{-1}\mathcal{Z})<\alpha^{-1}$, where $\mathds{A}(U) := A^\trp U A$ for $U\in\mathds{R}^{n\times n}$.
\end{enumerate}

If $\alpha\ge1$ and all eigenvalues of ${\alpha}A$ lay in the open unit disk then (i)--(v) are equivalent to $\alpha$-stochastic stability of $\Theta$.

If all eigenvalues of $A$ lay in the open unit disk then  (i)--(v) with $\alpha=1$  are equivalent to stochastic stability of $\Theta$.

\end{proposition}
\begin{proof} %
The equivalences among (i)--(v) are established in \cite[Th 3.3]{freiling}. We focus on the equivalence involving  $\alpha$-stability. In this proof, we refer to Lemma \ref{lemma:time.forward} applied to system $\Theta$, which amounts to set $B=0$ therein.

\medskip
[\emph{$\alpha$-stability/stability  $\Rightarrow$\/}].~ Suppose that $\Theta$ is $\alpha$-stable for $\alpha<1$. This implies necessarily that $\alpha^{\kappa} E_{x}[\|x(\kappa)\|^2]\to 0 $, and hence, for any $\Phi\in\mathds{R}^{n\times n}$ and $\theta\in\mathds{R}^n$ that
\begin{equation*}
\alpha^{\kappa} E_{x} \bigl[ \| x({\kappa}) \|^2_{\Phi} + \langle \theta, |x({\kappa})|  \rangle + \gamma \bigr]
\\
\le \alpha^{\kappa} E_{x} \bigl[ \| x({\kappa}) +\frac{1}{2}\Phi^{-1}\theta\sqcdot\mathcal{S}(x({\kappa}))| \|^2_{\Phi}+ \gamma \bigr]
\to 0
\end{equation*}
as $\kappa\to\infty$. 
From Lemma \ref{lemma:time.forward}, since $\kappa\to\mathcal{E}^{\alpha, \kappa}_2 (x(\cdotp))$ is monotone nondecreasing for sufficiently large $\kappa$,  
\begin{equation} \label{eq:bounded.energy}
    \lim_{\kappa\to\infty}    \mathcal{E}^{\alpha, \kappa}_{2,I} (x(\cdotp))\le  \lim_{\kappa\to\infty}\bigl(\| x \|^2_{P^{(\kappa)}_{0}} + \langle E_x [v_{0}^{(\kappa)}], x  \rangle + g^{(\kappa)}_{0} \bigr)\\
    \le \mathcal{E}^{\alpha}_2 (x(\cdotp))<\infty 
\end{equation}
holds for each $x(0) = x \in \mathds{R}^n$ with $P^{(\kappa)}(\cdotp)$, $v^{(\kappa)}(\cdotp)$ and $g^{(\kappa)}(\cdotp)$ and $Q=I$ in \eqref{algebraicequationsH2normcoef}, where we indicate the total number of stages by the superscript. From \eqref{eq:bounded.energy}, finiteness for each $x\in\mathds{R}^n$ and the fact that $\mathcal{L}^\alpha$ is linear-positive both imply that the solution $\kappa \to P_0^{(\kappa)}$ is monotone non-decreasing sequence for sufficiently large $\kappa$,  and converges to the solution of  \eqref{lyapunovkind.1} with $Q=I$. If this is true, (iii) is true, and  $(I-\mathcal{L}^\alpha)^{-1}(U)=Q\succ0$ has an uniquely defined solution $U\succ0$.

 When $\Theta$ is $\alpha$-stable for $\alpha \ge1$, there exist $c_0,c_1$ independent of $x(\cdotp)$ such that $\mathcal{E}^{\kappa,\alpha}_2(x(\cdotp)) \le c_0\|x-\xi_\kappa\|^2+ \kappa c_1\alpha^\kappa$ holds for bounded $\xi_\kappa$ and each $\kappa>0$. From the representation in Lemma \ref{lemma:time.forward} with $Q=I, B=0$ and $\Phi$, $\theta$  and $\gamma$ all set to zero, 
 \begin{equation}\label{eq:alpha.ge.1.bound}
     0\le \mathcal{E}^{\alpha, \kappa}_{2,I} (x(\cdotp)) = \| x \|^2_{P^{(\kappa)}_{0}} + \langle E_x [v_{0}^{(\kappa)}], x  \rangle + g^{(\kappa)}_{0}
    \\
    \le c_0\|x-\xi_\kappa\|^2+ \kappa c_1\alpha^\kappa
 \end{equation}
 
 Note from the assumption for $\alpha\ge1$ that
 $E_x [v_{0}^{(\kappa)}]$ is a bounded sequence, provided that there exists some upper bound matrix $L\succeq P^{(\kappa)}_{k}, 0\le k\le \kappa$. In this situation, one can evaluate,  
\begin{equation}\label{eq:uniform.bound.barv}
E_x [v_{0}^{(\kappa)}]= E_x \Bigl[\sum_{k=0}^{\kappa-1}\alpha^{k+1}\left( A^\trp\right)^k\mathcal{W}(P^{(\kappa)}_{k+1})\mathcal{S}(x_{k})\Bigr]\le \bar{v},\quad \forall \kappa\ge0
\end{equation}
where
\begin{equation}\label{eq:bound.exp.v0}
  \bar{v} := \alpha r_\sigma\left((I-\alpha A^\trp)^{-1}\right)\bigl|\mathcal{W}_d(L)\bigr|
\end{equation}

Now, by direct comparison, \eqref{eq:alpha.ge.1.bound} implies straightforwardly that 
\begin{gather}\label{eq:bound.P0.v0}
\| x \|^2_{P^{(\kappa)}_{0}} + \langle E_x [v_{0}^{(\kappa)}], x  \rangle  
=\|x -\zeta \|^2_{P^{(\kappa)}_{0}} -  \|\zeta\|^2_{P^{(\kappa)}_{0}}
\le c_0 \| x -\xi_\kappa\|^2
\\
g^{(\kappa)}_{0} =\sum_{k=0}^\kappa\alpha^k\tr(P^{(\kappa)}_k(\sigma\sigma^\trp + \sigma_x \sigma_x^\trp))
    \le\kappa \alpha^\kappa c_1 \label{eq:g_0(kappa)}
\end{gather}
with $\zeta = - (1/2)(P^{(\kappa)}_{0})^{-1} E_x [v_{0}^{(\kappa)}]$ for $P^{(\kappa)}_{0}$ invertible.  In the present setting, $\kappa\to P^{(\kappa)}_{0}$ defined by \eqref{recursivelyapunov} is such that $P^{(\kappa)}_{0}=\sum_{k=0}^{\kappa-1}\left(\mathcal{L}^\alpha\right)^k(I) \succeq I$, with $\left(\mathcal{L}^\alpha\right)^0 =I$, thus invertible.
The fact that $x$ is arbitrary, \eqref{eq:alpha.ge.1.bound}--\eqref{eq:g_0(kappa)}
tell us that  $P^{(\kappa)}(\cdotp)$ and $E_x [v_{0}^{(\kappa)}]$ are bounded sequences for every $\kappa>0$, and we can set $\xi_k=\zeta$. 

Since $\mathcal{L}^\alpha$ is linear-positive, it implies that the sequence $\kappa \to P_0^{(\kappa)}$ according with \eqref{recursivelyapunov} is monotone non-decreasing and converges as $\kappa\to\infty$ to the solution of  \eqref{lyapunovkind.1} with $Q=I$. If this is true, (iii) is true, and   $(I-\mathcal{L}^\alpha)^{-1}(U)=Q\succ0$ has an uniquely defined solution $L\succ0$. Set such a $L$ as the upper bound in need for \eqref{eq:uniform.bound.barv}.

For the C\`esaro's mean, one can apply \eqref{eq:alpha.ge.1.bound} with $\alpha$ set to one to get for the power norm that 
\begin{equation}\label{eq:stability.bound}
    \lim_{\kappa\to\infty} \frac{1}{\kappa}\mathcal{E}^{\alpha, \kappa}_{2,Q} (x(\cdotp)) \Big|_{\overset{\alpha =1}{Q=I}}=
    \\
    \lim_{\kappa\to\infty} \frac{1}{\kappa}\left(\| x \|^2_{P^{(\kappa)}_{0}} + \langle E_x [v_{0}^{(\kappa)}], |x|  \rangle + g^{(\kappa)}_{0}\right)
    < c_1
 \end{equation}
Since by assumption $r_\sigma(A)<1$ and reasoning as above with $\alpha=1$, one verifies that \eqref{eq:stability.bound} can only be satisfied for any $x\in\mathds{R}^n$ if 
\[
P^{(\kappa)}_0 = \sum_{k=0}^{\kappa-1}\mathcal{L}^k(I)\text{\quad and\quad  }E_x [v_{0}^{(\kappa)}]=E_x \Bigl[\sum_{k=0}^{\kappa-1}\left(A^\trp\right)^k\mathcal{W}(P^{(\kappa)}_{k+1})\mathcal{S}(x_{k})\Bigr]
\]
are bounded sequences, and thus,  
\begin{equation*}
      \lim_{\kappa\to\infty} \frac{1}{\kappa}\mathcal{E}^{\alpha, \kappa}_{2,Q} (x(\cdotp)) \Big|_{\overset{\alpha =1}{Q=I}}= \lim_{\kappa\to\infty} \frac{1}{\kappa}g^{(\kappa)}_{0}
      =  \lim_{\kappa\to\infty} \frac{1}{\kappa} \sum_{k=0}^\kappa\tr(P^{(\kappa)}_k(\sigma\sigma^\trp + \sigma_x \sigma_x^\trp))
   =\varpi(L) <c_1
\end{equation*}

The above implies that $\mathcal{L}$ is an inverse-positive operator and hence, the equivalences (i)--(v).

\medskip
[\emph{$\Rightarrow$  $\alpha$-stability/stability\/}].  
Suppose now that (iii) is satisfied, or equivalently, \eqref{lyapunovkind.1} is satisfied. Set $Q=I$ in \eqref{lyapunovkind.1} and denote $L= (I-\mathcal{L}^\alpha)^{-1}Q$.  Set  in Lemma \ref{lemma:time.forward} $P_\kappa=L$, $v_{\kappa} =0$ and $g_{\kappa} = 0$ with $B=0$. It yields, for any $x(0)=x$ and  $\alpha\ge 0$ that 
\begin{equation}
    \| x \|^2_{L} +  \langle E_{x}[v_0^{(\kappa)}], x  \rangle + g_{0}^{(\kappa)} =  \\
     \mathcal{E}^{\alpha, \kappa-1}_{2,I} (x(\cdotp)) + \alpha^{\kappa} E_{x} \bigl[ \| x({\kappa}) \|^2_{L} \bigr]
\label{eq:bounds.discounted}
\end{equation}
 where,
 \begin{equation*}
      v_0^{(\kappa)}=\sum_{k=0}^{\kappa-1} \alpha^{k+1}\left(A^\trp\right)^k\mathcal{W}(L)\mathcal{S}(x_k),
      \qquad
      g_{0}^{(\kappa)} = \sum_{k=0}^{\kappa-1} \alpha^{k+1}\varpi(L).
 \end{equation*}
Note that $v_0^{(\kappa)}$ is a well defined random vector for each $\kappa$, and we can evaluate, $ \lim_{\kappa\to\infty} |v_0^{(\kappa)}| 
    \le \bar{v}$, where $\bar{v}$ as in \eqref{eq:bound.exp.v0}.
That holds from (v), if $0\le\alpha\le1$, otherwise, from the assumption.

From \eqref{eq:bounds.discounted} and the fact that here $L$ is invertible, we get that 
 \begin{multline}\label{eq:bound.for.stability}
      0\le \mathcal{E}^{\alpha, \kappa-1}_{2,I} (x(\cdotp))\le \| x \|^2_{L} +  \langle E_{x}[v_0^{(\kappa)}], x  \rangle + g_{0}^{(\kappa)}
      \\
      \le 
      \| x -\xi_\kappa \|^2_{L} +  \frac{\alpha^{\kappa+1} -\alpha}{\alpha-1}\varpi(L)
      \le \bigl\| x -\xi_\kappa\|^2_{L} +
      \alpha^\kappa \frac{\alpha\varpi(L)}{\alpha-1}
\end{multline}
holds, with $\xi_\kappa= -\frac{1}{2}L^{-1}E_{x}[v_0^{(\kappa)}]<\infty$.
Hence, we conclude that $\Theta$ is $\alpha$- stochastically stable for $\alpha\ge0$, since \eqref{eq:bound.for.stability} provides that
\begin{align*}\allowdisplaybreaks
     &\text{For $\alpha<1$, }\mathcal{E}^{\alpha}_{2} (x(\cdotp))\le\infty, 
     \\
     &\text{For $\alpha>1$, } c_0 =\lambda^+(L), \xi_\kappa = -\frac{1}{2}L^{-1}E_{x}[v_0^{(\kappa)}], \text{ and $c_1= \frac{\alpha\varpi(L)}{\alpha-1}$,}
     \\
     &\text{For $\alpha=1$, $c_0$, $\xi_\kappa$ as above,  and by identifying $g_0^{(\kappa)}$, $c_1= \varpi(L)$}.
\end{align*}

 For the C\`esaro mean it readily follows that
\begin{equation*}
    \limsup_{\kappa\to\infty}\frac{1}{\kappa}\mathcal{E}^{1,\kappa}_{2,I} (x(\cdotp)) \le \lim_{\kappa\to\infty}\frac{1}{\kappa}\left(\bigl\| x + \frac{1}{2}L^{-1}\bar{v}\bigr\|_L^2+  \kappa\varpi(L)\right)=\varpi(L)
\end{equation*}
Now write $\mathcal{E}^{1,\kappa-1}_{2,I}(x(\cdotp))+E[\|x(\kappa)]\|^2_L]= \mathcal{E}^{1,\kappa}_{2,I}(x(\cdotp))+E[\|x(\kappa)]\|^2_{\Delta L}]$ with $\Delta L = L-\mathcal{L}(L)=I$,  and from \eqref{eq:bounds.discounted}, it follows that 
\begin{multline*}
    \liminf_{\kappa\to\infty}\frac{1}{\kappa}\Bigl(\mathcal{E}^{1,\kappa}_{2,I}(x(\cdotp))+E[\|x(\kappa)]\|^2]\Bigr)
    =\liminf_{\kappa\to\infty}\frac{1}{\kappa}\bigl(\| x \|^2_{L} +  \langle E_{x}[v_0^{(\kappa)}], x  \rangle + g_{0}^{(\kappa)}\bigr)
    \\
    \geq -\lim_{\kappa\to\infty}\frac{1}{\kappa}\left(\bigl\| x + \frac{1}{2}L^{-1}\bar{v}\bigr\|_L^2\right) + \varpi(L)=\varpi(L)
\end{multline*}
which shows that $ \mathcal{P}_{2} = \varpi(L)$ and $\Theta$ is stochastically stable, completing the proof.
\end{proof}

From the proof of Proposition \ref{prop:freiling.stab}, we have the following corollary.

\begin{corollary} \label{coro:stability.discounted}
System $\Theta$ in (\ref{CSVIUsystem_no_control}) is $\alpha$-stochastically stable iff
\begin{enumerate}
    \item[i)] $\alpha<1$, $L\succ0$ is the solution of \eqref{lyapunovkind.1} for  $Q\succ0$,
    \item[ii)]
     $\alpha\ge1$, $L\succ0$ is the solution of \eqref{lyapunovkind.1} for  $Q\succ0$ provided that all eigenvalues of ${\alpha}A$ lies in the open unit disk.
\end{enumerate} 
System $\Theta$ in (\ref{CSVIUsystem_no_control}) is stochastically stable iff (ii) holds for $\alpha=1$.

Moreover, when $0\le\alpha<1$ its Q-mean energy is expressed by 
\begin{equation}    \label{eq:limite_discounted}
      \mathcal{E}^\alpha_{2,Q}(x(\cdot)) = \lim_{\kappa \to \infty}  E \Big[\sum_{k=0}^{\kappa} \alpha^k\| x(k) \|_Q^2 \Big] 
      \\
      = \| x \|^2_{L} +  \langle E_{x}[v], x  \rangle + g_{0} 
\end{equation}
where, 
 $
      v =\lim_{\kappa\to\infty}\sum_{k=0}^\kappa \alpha^{k+1}\left(A^\trp\right)^k\mathcal{W}(L)\mathcal{S}(x_k)
$ 
 and $g_0 = {\alpha}\varpi(L)/{(1-\alpha)} $.
 
The C\`esaro mean for the $Q$-mean power is given by
\begin{equation}    \label{eq:limite_long.run}
      \mathcal{P}_{2,Q}(x(\cdot)) = \lim_{\kappa \to \infty} \frac{1}{\kappa} E \Big[\sum_{k=0}^{\kappa} \| x(k) \|_Q^2 \Big] 
      = \varpi(L)
\end{equation}

\end{corollary}

\begin{proof}
For some $Q\succ0$ in \eqref{lyapunovkind.1}, as in \eqref{eq:bounds.discounted}, set in Lemma \ref{lemma:time.forward} with $B=0$, $P_\kappa=L$, $v_{\kappa} = E[v|x_\kappa]$ and $g_{\kappa} = 0$ 
to get,
\begin{equation*} 
 \mathcal{E}^{\kappa,\alpha}_{2,Q}(x(\cdot)) = 
\| x \|^2_{L}    + \langle E_{x}[v], x\rangle  + \frac{\alpha^{\kappa-1}-1}{\alpha - 1}\alpha\varpi(L)
 \end{equation*}
When $0\le\alpha<1$, take the limit as $\kappa\to\infty$ to get \eqref{eq:limite_discounted}. 

For the $Q$-mean power, we rewrite \eqref{eq:bounds.discounted} with $\alpha=1$ to get, similarly to the proof of Proposition \ref{prop:freiling.stab}, 
\begin{equation*}
    \lim\sup_{\kappa\to\infty}\frac{1}{\kappa}\mathcal{E}^{ \alpha, \kappa}_{2,Q} (x(\cdotp))\Big|_{\alpha =1} \le \varpi(L) \le 
    \lim\inf_{\kappa\to\infty}\frac{1}{\kappa}\mathcal{E}^{ \alpha, \kappa}_{2,Q} (x(\cdotp))\Big|_{\alpha =1}
\end{equation*}
\end{proof}

\subsection{Stochastic Detectability}\label{sec:detectability}
This short but too important section allows us to weaken the requirement that $Q\succ0$ regarding the assurance of stochastic stability, to some qualified matrix $Q\succeq0$. 
Here we introduce an important notion that can guarantee  $\alpha$-stability when the perturbed Lyapunov matrix equation \eqref{lyapunovkind.1} possesses a positive semidefinite solution, specially that associate with $Q=C^\trp C$. 

\begin{definition}\label{def:detability}
System $\Theta$ is $(C,\alpha^{1/2}A,\alpha\mathcal{Z})$ $\alpha$-detectable if  there is a matrix $G \in \mathds{R}^{n \times p}$ such that $\alpha^{1/2}(A+GC)$ is $d$-stable relative to $\alpha\mathcal{Z}$. 

We often call this notion $(C,\mathcal{L}^\alpha)$-detectable for some $\alpha\ge0$ or simply $(C,\mathcal{L})$-detectable when $\alpha=1$.
\end{definition}

In view of the equivalences in Proposition \ref{prop:freiling.stab}, the next result can be announced.

\begin{theorem}\label{thm:detectability}
Suppose that for some $\alpha\ge0$,
\[
(I-\mathcal{L}^\alpha)(U) =C^\trp C 
\]
has a solution $L\succeq0$. If $\Theta$ is $(C,\mathcal{L}^\alpha)$-detectable and if either,
\begin{enumerate}
    \item[i)] $\alpha<1$, or,
    \item[ii)] $\alpha\ge1$ and all eigenvalues of $\alpha A$ lies in the open unit disk,
\end{enumerate}
then  $\Theta$ is $\alpha$-stochastically stable and conversely. 
 \begin{enumerate}
     \item[iii)] When (ii) holds for $\alpha=1$, $\Theta$ is stochastically stable and conversely. 
 \end{enumerate}

\end{theorem}

\begin{proof}
        (i) It is a straighforward conclusion of Lemma\,3.7 of \cite{freiling} and Proposition \ref{prop:freiling.stab}. The condition (ii) adds the requirement of Proposition \ref{prop:freiling.stab}. The last assertion is a direct consequence of 
        (ii).
\end{proof}

\begin{corollary}[Detectability and Stochastic Stability]\label{coro:detectability.stability}
Suppose that system $\Theta$ is $(C,\mathcal{L}^\alpha)$-detectable. Then,
\begin{enumerate}
    \item[i)] 
    System $\Theta$ is $\alpha$-stochastically stable  iff \eqref{lyapunovkind.1} with $Q=C^\trp C$ has a solution $L\succeq0$, and if  $\alpha\ge1$, condition (ii) in Theorem \ref{thm:detectability} should also hold.
    
     \item[ii)] 
     If $\lim_{\kappa\to\infty}\mathcal{E}_{2,Q}^{\alpha,\kappa}(x(\cdot))<\infty$ for $Q=C^\trp C$ and any $x(0)$, then $\Theta$ is $\alpha$-stable. 
    \item[iii)] If system $\Theta$ is $\alpha$-stochastically stable its Q-mean energy is expressed by eqn. \eqref{eq:limite_discounted}, where $g_0$ may be unbounded for $\alpha\ge1$.  
    \item[iv)] If system $\Theta$ is stochastically stable, its Q-mean power appears in eqn. \eqref{eq:limite_long.run}. 
 \end{enumerate}
\end{corollary}

\begin{remark}
Corollary \ref{coro:detectability.stability} connects finiteness of  a $Q$-energy measurements of type $\mathcal{E}_{2,Q}^\alpha(x(\cdotp))$ with the corresponding $\alpha$-stochastic stability notion of system $\Theta$, whenever $\Theta$ is $(C,\mathcal{L}^\alpha)$-detectable. When $\alpha\ge1$ the $Q$-measurement may be unbounded; however, $\mathcal{E}_{2,Q}^{\alpha,\kappa}(x(\cdotp))$ is bounded for each $\kappa$ according to Definition \ref{stability_definition1} (ii).  

In addition, stability and finiteness are tied to the solution of \eqref{lyapunovkind.1}, and provided that $\alpha<1$, it comes as an exact parallel to deterministic linear system \cite{Kailath} as well as some stochastic problems \cite{Davis,EFCosta,bib:2,bib:11}.  
However, when $\alpha\ge1$, a further requirement is in need, the one in Theorem \ref{thm:detectability} (ii).
\end{remark}

\section{Abel and  C\'esaro mean and the Long-Run Approach}
\label{subsec:long-run.ana}

To delve with system behavior in the power norm measurement in \eqref{energymeasure2} or in the energy measurement  \eqref{energymeasure1} when $\alpha>1$ is the main aim of this section.   
 When $\alpha=1$, we refer generically to the long-run approach, and the relation between the Abel and C\`esaro means are of interest here.
 
 The energy measurement, when coupled with a discount factor ($\alpha<1$), leads to a stability sense which implies that $\alpha^{k}\|x_k\|^2 \to 0$ in the mean-square sense, c.f. \cite{bib:4,bib:3}. This might be too feeble, e.g., for the sake of comparison, it would imply only that all eigenvalues of $\sqrt{\alpha}A$ lies in the open unit disk, not even guaranteeing the stability of the deterministic LTI system if the noise would subside. 
 
 In comparison, when $\alpha\ge1$ and large, the steady state behavior has an increasing impact on the $Q$-energy measurement as the time horizon increases, magnifying requirements on the stability and restraints on drifts away from the origin neighborhood. Such analysis encompasses the system's performance in the far future (or on the steady state), bringing substantial requirements on the asymptotic of state and output processes. 
 
In this section we show the state convergence to a known value in the mean-square sense, and provide an exponential bound on convergence error in this same sense. These stands as clear improvements on the convergence notion so far obtained for the discounted case, see Lemma \ref{lemm:convergence.alpha.ge1}.

To connect the Abel's and  C\'esaro's mean, we choose the ``vanishing discount/ counter-discount" approach. The idea is to treat the long-run average power norm as the limit of a sequence of discounted and counter-discounted energy norm problems, as the discount factor  $\alpha \to 1$, and rely on convergence, see \cite{bib:8,bib:9} for the standard vanishing discount approach.

Suppose that \eqref{lyapunovkind.1} has a solution for some $\alpha\ge0$ and when $\alpha>1$, there holds $r_\sigma(A)<\alpha^{-1}$. It is a simple matter to conclude from Proposition \ref{prop:freiling.stab} that \eqref{lyapunovkind.1} will have a solution for any $\alpha'$ such that $\alpha'\le \alpha$, and we denote  $\bar{\alpha}: =\sup \{\alpha :$ \eqref{lyapunovkind.1} has a solution  and  $r_\sigma(A)<\alpha^{-1}\}$.

Note also from the definition that if a pair $(C,\mathcal{L}^\alpha)$ is detectable for some $\alpha$, the pair $(C,\mathcal{L}^{\alpha'})$ is necessarily detectable when $\alpha'\le \alpha$.

 \begin{lemma} 
 Set $\mathcal{L}\equiv \mathcal{L}^{\alpha=1}$ 
 and suppose that there exists a solution $L^\alpha$ to the equation
 \begin{equation}\label{convergent_lyapunov}
 (I-\mathcal{L})(U)=Q
\end{equation}
 Let $\alpha_n \to 1, n \geq 0$ be a  sequence, such that the corresponding sequence of solutions of
 \begin{equation}\label{eq:sequence_Lyapunov}
     (I-\mathcal{L}^{\alpha_n}) (U) = Q 
 \end{equation}
 are $L^{\alpha_n}\succeq0$ for each $n$.
Then  $L^{\alpha_n}\to L$ in the semi positive definite sense for $L \succeq 0$, the unique solution of \eqref{convergent_lyapunov}.
\label{lema_emergent_lyapunovs}
\end{lemma}

\begin{proof} 
First we show that the sequence  $L^{\alpha_n}\succeq0$ is non-void. This is an easy task, since if \eqref{convergent_lyapunov} is satisfied, there exists some  $\tilde{\alpha}>1$ such that $(I-\mathcal{L}^{\alpha}) (U) = Q$ has a solution $L^{\alpha}\succeq0$ for any $\alpha<\tilde{\alpha}$. Due to the fact that $\mathcal{L}^\alpha$ is made up by the sum of two  linear-positive and monotone operators, it follows that $\mathcal{L}^{\alpha_n}$ is an inverse-positive operator according with Proposition \ref{prop:freiling.stab} for each $n\ge0$, and the solutions $L^{\alpha_n} =(I-\mathcal{L}^{\alpha_n})^{-1}Q\succeq0$ are unique.  
Note also that if $0\le\alpha_1<\alpha_2$,
\begin{equation}\label{eq:monotone}
    \mathcal{L}^{\alpha_2}(U)\succeq \mathcal{L}^{\alpha_1}(U).
\end{equation}
Subtracting successive equations \eqref{eq:sequence_Lyapunov} and using their corresponding solutions, one gets that
\[
(\mathcal{L}^{\alpha_n}-I)(L^{\alpha_n}-L^{\alpha_{n-1}})+\mathcal{L}^{\alpha_n}(L^{\alpha_{n-1}}) - \mathcal{L}^{\alpha_{n-1}}(L^{\alpha_{n-1}})=0
\]

Now, we highlight two subsequences of $\{\alpha_n\}$ such that
$\alpha^-_m= \min\{\alpha_n: n\ge m\}$ and $\alpha^+_m= \max\{\alpha_n: n\ge m\}, m\ge0$. These are monotone sequences, and we reduce the analysis to either monotone increasing $\alpha^-_n\uparrow1$ or decreasing  $\alpha^+_n\downarrow1$ sequences, without loss.

For the increasing sequence it implies from \eqref{eq:monotone}  that $\mathcal{L}^{\alpha^-_{n}}(U)- \mathcal{L}^{\alpha^-_{n-1}}(U)\succeq0, \forall U\in\mathds{S}^{n+}, n\ge0$. This, in turn, provides that
$L^{\alpha^-_n}-L^{\alpha^-_{n-1}}\succeq0$. Hence, $0\preceq L^{\alpha_0}\preceq\ldots \preceq L^{\alpha_n}\uparrow L^1$, where $L^1$ satisfies \eqref{convergent_lyapunov} and from uniqueness, $L^1\equiv L$, which  shows the result for the monotone increasing subsequence. 

For the monotone decreasing subsequence, the ordering $\alpha^+_{n-1}\ge \alpha^+_n\ge \ldots$ is reversed, and also are all comparisons in the positive semi-definite sense developed above, completing the proof.
\end{proof}

\begin{lemma} \label{lemma:converegence.alfa.qqer}
 Suppose that  $\Theta$ in (\ref{CSVIUsystem_no_control}) is $(Q^{1/2} \mathcal{L}^\alpha)$-detectable, there exists a solution  $L\succeq0$ for  \eqref{lyapunovkind.1} and there holds $r_\sigma(A)<\alpha^{-1}$ for some $0<\alpha<\bar{\alpha}$. Then for each $\kappa\ge0$,
\begin{equation}\label{eq:finitiness.alpha.ge1}
\Bigl| \sum_{k=0}^{\kappa} \alpha^{k} E\bigl[ \|x_k\|^2_Q-\alpha\varpi(L) \bigr] \Bigr| \le \|x_0\|_L^2+\langle\bar{v},|x_0|\rangle   
\end{equation}
for some $\bar{v}\in\mathds{R}^n$. Moreover,
\begin{equation}\label{eq:bound.norm.state}
    E[\|x_\kappa-\xi_\kappa \|^2]\le c_0(1+2\alpha^{-\kappa}) + c_1\alpha^{-\kappa}\|x_0 \|^2 
\end{equation}
for some $c_0,c_1>0$ and bounded $\xi_\kappa\in\mathds{R}^n, \forall\kappa$.

\end{lemma}

\begin{proof}
 Let us consider $P^{(\kappa)}_k$, $v^{(\kappa)}_k$  and  $g^{(\kappa)}_k$,  for  $k=0,1,\ldots,\kappa$,  the solution of eqns. \eqref{algebraicequationsH2normcoef} in Lemma \ref{lemma:time.forward} with $B=0$. The superscripts refer to horizon number $\kappa$, and we also set $\Phi=P^{(\kappa)}_\kappa=0$,  $\theta=v^{(\kappa)}_\kappa=0$ and $\gamma=g^{(\kappa)}_\kappa=0$.
 In this situation, the lemma yields the representation 
\begin{equation} \label{eq:sum_and_multiply_by_(-1)}
      E \big[\sum_{k=0}^{\kappa-1} \alpha^k \|x_k\|^2_Q\big] = 
      \| x \|^2_{P^{(\kappa)}_{0}} +  \langle E_{x}[v^{(\kappa)}_0], x  \rangle + g^{(\kappa)}_{0}, \quad x_0=x,
\end{equation}
with $g^{(\kappa)}_{0} = \sum_{k=0}^{\kappa-1} \alpha^{k+1} \varpi(P^{(\kappa)}_k)$.

Since system $\Theta$  is $(Q^{1/2},\mathcal{L}^\alpha$)-detectable and $r_\sigma(A)<\alpha^{-1}$, Theorem \ref{thm:detectability} implies that $\Theta$ is $\alpha$-stochastically stable in the sense of Definition \ref{stability_definition1}. Moreover, by construction note that 
 $0\preceq P^{(\kappa)}_k\preceq P^{(\kappa+1)}_k\preceq \cdots\preceq P^{(\infty)}_k = L$, the solution of \eqref{lyapunovkind.1} holds for each $k\leq\kappa$.
Also, 
\begin{equation*}
    v_k^{(\kappa)} = \sum_{n=0}^{\kappa-k} \alpha^{n+1}\left(A^\trp\right)^n \mathcal{W}(P^{(\kappa)}_{n+k})\mathcal{S}(x_{n+k}) 
\end{equation*} 
and as in the proof of Proposition \ref{prop:freiling.stab},
\begin{equation}\label{eq:bar.v}
|v_k^{(\kappa)}|\le  |v_k^{(\infty)}|\le \bar{v} = \alpha r_\sigma\left((I-\alpha A^\trp)^{-1}\right)\bigl|\mathcal{W}_d(L)\bigr|
\end{equation}

Hence, from \eqref{eq:sum_and_multiply_by_(-1)}
\begin{multline*}
     \lim_{\kappa\to\infty}\Bigl|\sum_{k=0}^{\kappa-1} \alpha^k \bigl( E [\|x_k\|^2_Q] - \alpha \varpi(P^{(\kappa)}_k)   \bigr) \Bigr|
     =\Bigl|\sum_{k=0}^{\infty} \alpha^k \bigl(   E [\|x_k\|^2_Q] -\alpha\varpi(L)  \bigr) \Bigr|
     \\
     =\lim_{\kappa\to\infty}\Bigl| \| x_0 \|^2_{P^{(\kappa)}_{0}} +  \langle E_{x}[v^{(\kappa)}_0], x_0  \rangle \Bigr|
     \le \|x_0\|^2_{L} + \langle \bar{v},|x_0|\rangle
\end{multline*}
which shows \eqref{eq:finitiness.alpha.ge1}.

 For the bound on the second moment in \eqref{eq:bound.norm.state}, let us denote for  $k=0,1,\ldots,\kappa$,
\begin{equation}\label{eq:var.V.alpha}
 \Psi_k(x_k) := \| x_k \|^2_{P^{(\kappa)}_k} +  \langle E [v^{(\kappa)}_k|x_k], x_k  \rangle, 
\end{equation}
 where, similarly to the above, $P^{(\kappa)}_k$ and  $v^{(\kappa)}_k$ are the sequences produced by \eqref{recursivelyapunov} and \eqref{recursivevector.2} respectively, in Lemma \ref{lemma:time.forward} with $B=0$. But here, one sets $\Phi=P^{(\kappa)}_\kappa=L$ and for each  $k=0,\ldots,\kappa$,
\begin{equation*}
 v^{(\kappa)}_k= 
    \sum_{n=0}^{\infty} \alpha^{n+1}\left(A^\trp\right)^n \mathcal{W}(L)\mathcal{S}(x_{n+k}) 
\end{equation*} 
Clearly, if $\theta = v^{(\kappa)}_\kappa$ as above, it sets $v_k = v^{(\kappa)}_k$ as above, for each $k=0,\ldots,\kappa$. Write $\Psi^\alpha_k = \Psi^\alpha(x_k)$ for short, and having in mind Lemma \ref{lemma:time.forward}, and that $k\to x_k$ is a Markovian process, we write,
\begin{multline}\label{eq:sum.varV.1}
\allowdisplaybreaks
        \alpha^\kappa E[\Psi_\kappa -\Psi_0]=E \Big[ \sum_{k=0}^{\kappa-1} \alpha^k (\alpha \Psi_{k+1} - \Psi_{k}) \Big]= \\
         \quad E \Big[ \sum_{k=0}^{\kappa-1} \alpha^k E\Big( x^\trp ( \mathcal{L}^\alpha (L) - L )x+ \alpha \varpi(L)+\\
         \langle \alpha A^\trp v_{k+1}  + \alpha \mathcal{W}(L)\mathcal{S}(x) - v_{k} ,x \rangle \big| x_k=x \Big) \Big]=
        E \Big[ \sum_{k=0}^{\kappa-1} \alpha^k \big(  \alpha \varpi(L) - \|x_k\|^2_Q\big) \Big]
\end{multline}
On other hand,
\begin{equation}\label{eq:sum.varV.2}
        \alpha^\kappa E[\Psi_\kappa -\Psi_0]=
        \alpha^\kappa E \bigr[ \|x_\kappa\|^2_{L}  +  \langle  v_\kappa, x_\kappa \rangle\bigl] -\|x_0\|^2_{L} - \langle E [v_0|x_0], x_0  \rangle
\end{equation}

Under the assumptions, one can  eventually set $Q\succ0$, which implies that the solution  of  \eqref{lyapunovkind.1}, $L=\sum_{n=0}^\infty (\mathcal{L}^\alpha)^n(Q)\succ 0$ (with $(\mathcal{L}^\alpha)^0=I$).
Thus, set $\xi_0=-\frac{1}{2}(L)^{-1}E[v_0|x_0]$ and $\bar{\xi} = -\frac{1}{2}(L)^{-1}\bar{v}$ with $\bar{v}$ as in \eqref{eq:finitiness.alpha.ge1}. From \eqref{eq:sum.varV.1} and \eqref{eq:sum.varV.2}, we can apply \eqref{eq:finitiness.alpha.ge1} to get,
\begin{multline*}
    \alpha^\kappa E \bigr[ \|x_\kappa\|^2_{L}  +  \langle  v_\kappa, x_\kappa \rangle\bigl] = 
    \\
    =\|x_0\|^2_{L} + \langle E [v_0|x_0], x_0  \rangle + 
    E \Big[ \sum_{k=0}^{\kappa-1} \alpha^k \big(  \alpha \varpi(L) - \|x_k\|^2_Q\big) \Big]
    \\
    \le  \|x_0-\xi_0\|^2_{L} + \|x_0 - \bar{\xi}\sqcdot\mathcal{S}(x_0)\|_{{L}}^{2}
\end{multline*}
In addition, note that
\begin{equation*}
    E \bigr[ \|x_\kappa\|^2_{L}  +  \langle  v_\kappa, x_\kappa \rangle\bigl] \ge
    E \bigr[ \|x_\kappa\|^2_{L}  -  \langle  \bar{v}, |x_\kappa| \rangle\bigl] \ge
    E \bigr[ \bigl\|x_\kappa -\bar{\xi}\sqcdot\mathcal{S}(x_\kappa)\bigr\|^2_{L}  -  \|\bar{\xi}\|^2_{L}\bigl]
\end{equation*}
 Therefore, 
\begin{multline*}
     E \bigr[ \|x_\kappa  -\bar{\xi}\sqcdot\mathcal{S}(x_\kappa)\|^2_{L} \bigl] \le   \|\bar{\xi}\|^2_{L} + \alpha^{-\kappa}\bigl(\|x_0-\xi_0\|^2_{L} +  \|x_0 - \bar{\xi}\sqcdot\mathcal{S}(x_0)\|_{{L}}^{2}\bigr)
    \\
    \le (1+2\alpha^{-\kappa})\|\bar{\xi}\|^2_{L} + 2\alpha^{-\kappa}\|x_0\|^2_{L}
\end{multline*}
Finally, denote $\bar{\xi}_\kappa=\bar{\xi}\sqcdot\mathcal{S}(x_\kappa)$ to get
\begin{multline*}
    E \bigr[ \|x_\kappa -\bar{\xi}_\kappa\|^2 \bigl] \le
    \\
      \bigr((1+2\alpha^{-\kappa})\|\bar{\xi}\|^2_{L} 
    + 2\alpha^{-\kappa}\|x_0\|^2_{L}\bigl)/\lambda^-(L)
    \le  c_0(1+2\alpha^{-\kappa}) + c_1\alpha^{-\kappa}\|x_0 \|^2
\end{multline*}
 Since $\|\bar{\xi}_\kappa\|\le \|\bar{\xi}\|<\infty,\forall \kappa$, the bound in \eqref{eq:bound.norm.state} is demonstrated.
\end{proof}

The evaluations in Lemma \ref{lemma:converegence.alfa.qqer} imply the following bound for the second moment and convergence notion.

\begin{lemma}\label{lemm:convergence.alpha.ge1} Suppose that  $\Theta$ in (\ref{CSVIUsystem_no_control}) is $(Q^{1/2} \mathcal{L}^\alpha)$-detectable, there exists a solution  $L\succeq0$ for  \eqref{lyapunovkind.1} and there holds $r_\sigma(A)<\alpha^{-1}$ for some $1\le \alpha <\bar{\alpha}$. Then, \begin{enumerate}
    \item[i)] 
    $E[\|x_k\|^2], k\ge0$ is bounded,
    \item[ii)] 
    $x_k$ converges in the m.s.-sense, namely, $E\bigl[ \|x_k\|^2_Q \bigr]\to \alpha\varpi(L)$ as $k\to\infty$,
    \item[iii)]
    For any $0<\alpha<\bar{\alpha}$, 
    \begin{equation}
        \left|E\bigl[ \|x_k\|^2_Q-\alpha\varpi(L) \bigr] \right|\le 2\alpha^{-k}\left( \|x_0\|_L^2+\langle\bar{v},|x_0|\rangle\right)
    \end{equation}
    holds for all $k>0$ and some $\bar{v}\in\mathds{R}^n$.
\end{enumerate} 

 \end{lemma}
 
\begin{proof} (i) From \eqref{eq:bound.norm.state} with $\alpha\ge1$ if follows that $ E[\|x_\kappa-\xi_\kappa \|^2]$ is bounded for each $\kappa$, and hence, $ E[\|x_\kappa \|^2]$ is also bounded.
Let us denote $s_N = \sum_{k=0}^{N} \alpha^{k}E\bigl[ \|x_k\|^2_Q-\alpha\varpi(L) \bigr]$ for $0<\alpha<\bar{\alpha}$. From \eqref{eq:finitiness.alpha.ge1}, one has that $|s_N|\le c(x_0) := \|x_0\|_L^2+\langle\bar{v},|x_0|\rangle $, valid for all $N$. Subtracting two values we get that $|s_{N_1} - s_{N_2}|\le 2c(x_0)$ and in particular, by setting $N_1 =N$ and $N_2=N-1$ we get that
\[
-2c(x_0)\le \alpha^N E\left[ \|x_N\|^2_Q-\alpha\varpi(L) \right]\le 2c(x_0)
\]
which shows (iii).

The above also implies for $1<\alpha <\bar{\alpha}$ that $E\bigl[ \|x_k\|^2_Q \bigr]\to \alpha\varpi(L)$ as $k\to\infty$. Since eventually, we can choose $Q=I$, this shows state convergence in the mean-square sense to the corresponding value $\alpha\varpi(L)$, with $L=(I-\mathcal{L}^\alpha)^{-1}(I)$.

Assertion (ii) then follows from (iii) straightforwardly when $\alpha>1$. When $\alpha=1$, Corollary \ref{coro:stability.discounted} and particularly \eqref{eq:limite_long.run}, express the validity of (ii), which completes the proof.
\end{proof}

\begin{remark}[Exponential convergence]\label{rem:convergence.properties}
When the energy problem can be solved for $\alpha>1$, Lemma \ref{lemm:convergence.alpha.ge1} adds to the subtleties on the long-run state behavior regarding mean-square finiteness and convergence to a known value.  Item (iii) gives a geometric decay measure of the m.s.-error in such a way that, 
\[
 \bigl| E\bigl[ \|x_\kappa\|^2_Q - \alpha\varpi(L)\bigl] \bigr| \le \beta \alpha^{-\kappa} 
\]
for some $\beta$ depending on $x_0$.
\end{remark}

Lemmas \ref{lemma:converegence.alfa.qqer} and \ref{lemm:convergence.alpha.ge1}, and the next proposition pave the way to another important result of this section.

\begin{proposition}[Tauberian Theorems, \cite{bib:6}]\label{prop.tauberian} Let  $\{ a_n \}$ be a  sequence of real numbers and $s_N(\alpha)=\sum_{n=0}^N\alpha^n a_n$. Set  $f_N(\alpha)=(1-\alpha)s_N(\alpha)$ and $f(\alpha)=\lim_{N\to\infty}f_N(\alpha)$. Denote also $s_N=s_N(1)$ and $h_N= s_N/(N+1)$.
\begin{enumerate}

    \item[i)]\emph{Hardy and Littlewood Theorem.\/} Suppose that $\{ a_n \}$ is a bounded sequence and let $\lim _{\alpha \rightarrow 1}f(\alpha)=a$. 
    Then, $\lim_{N \rightarrow \infty} h_{N}=a$.
    \item[ii)]
Let $\{ a_n \}$ be a non negative sequence and $\alpha \in (0,1)$. Then,
\begin{equation*}\allowdisplaybreaks
\liminf_{N \to \infty} h_N \leq \liminf_{\alpha \uparrow 1} f(\alpha)  \leq 
\\
\limsup_{\alpha \uparrow1} f(\alpha) \leq \limsup_{N \to \infty} h_N
\end{equation*}
\end{enumerate}
\label{tauberiantheorem}
\end{proposition}

Next, a vanishing discount/counter-discount result is provided that links the two types of energy and power measurements studied here. 
Note that when \eqref{convergent_lyapunov} has a solution, there exists $\bar{\alpha}$ as in Lemma \ref{lemma:converegence.alfa.qqer} with $\bar{\alpha}>1$.

\begin{theorem}(Stochastic Stability: Abel and C\`esaro means)  \label{thm:analysis:longrunnorm2}
Suppose that $\Theta$ in (\ref{CSVIUsystem_no_control}) is $(Q^{1/2}, \mathcal{L})$-detectable and $r_\sigma(A)<1$. Then $L \succeq0$ is the solution of \eqref{convergent_lyapunov} iff  system $\Theta$ is stochastically stable and its Q-mean average power is given by \eqref{eq:limite_long.run}.
Moreover, for $\alpha<\bar{\alpha}$ with $\bar{\alpha} = \sup\{ \alpha: L^\alpha=(I-\mathcal{L}^\alpha)^{-1}(Q)\succeq0$ and $r_\sigma(A)<\alpha^{-1}\}$,
\begin{multline}\label{eq:vanishing.general}\allowdisplaybreaks
       \lim_{\alpha \to 1} (1 - \alpha)
       \sum_{k=0}^\infty \alpha^{k}\left( E[\|x_k\|^2_Q]-\alpha\varpi(L^\alpha)\right)
    \\
    =\lim_{\kappa\to\infty}\frac{1}{\kappa}\sum_{k=0}^{\kappa-1} \left( E[\|x_k\|^2_Q]-\varpi(L)\right) 
    =\mathcal{P}_{2,Q}(x(\cdot)) -\varpi(L)=0
 \end{multline}
\end{theorem}

\begin{proof}
The first assertion is a restatement of Theorem \ref{thm:detectability} (iii) for stochastically stability and the Q-mean average power of $\Theta$.

Set $a_k :=  E[\|x_k\|^2_Q]-\alpha\varpi(L^\alpha)$, and for any sequence $\alpha_n\to 1$, let us divide it in two subsequences, $\{\alpha^+_n\}=\{\alpha_m:1\le \alpha_m < \bar{\alpha}\}$ and $\{\alpha^-_n\}=\{\alpha_m:0\le \alpha_m<1\}$. 
Note from the assumptions that $\bar{\alpha}>1$. To prove the result, we show that both sequences set \eqref{eq:vanishing.general} true.

First, for sequence $\alpha^+_n\downarrow1$, let us consider 
 that there exists $L^{\alpha_n^+} \succeq0$ and $r_\sigma(A)<(\alpha_n^+)^{-1},\forall n\ge0$. 
 Lemma \ref{lemm:convergence.alpha.ge1} tell us 
  that the sequence $\{E[\|x_k\|^2_Q]\}$ is bounded when $\alpha \ge1$, and hence, $\{a_k\}$ is bounded, since $\varpi(L^{\alpha^+_n})$ is obviously bounded. We can apply Proposition \ref{prop.tauberian} (i) to sequence $\{\alpha^+_n\}$ to conclude 
that \eqref{eq:vanishing.general} holds for any $\alpha_n\downarrow1$.

Regarding sequence $\{\alpha^-_n\}$, one has from \eqref{eq:finitiness.alpha.ge1}  that 
\begin{multline*}
\lim_{\kappa\to \infty} \Bigl| \sum_{k=0}^{\kappa} E\Bigl[\alpha^{k}\left( \|x_k\|^2_Q-\alpha\varpi(L^\alpha) \right) \Big] \Bigr| 
= \Bigl| \mathcal{E}^{\alpha}_{2,Q}(x(\cdot))-\frac{\alpha}{1-\alpha}\varpi(L^\alpha) \Bigr|
\\
\le \|x_0\|_{{L}^{\alpha}}^{2}+\langle\bar{v},|x_0|\rangle  
\end{multline*}
or, 
\[
\mathcal{E}^{\alpha}_{2,Q}(x(\cdot)) \le \|x_0\|_{{L}^{\alpha}}^{2}+\langle\bar{v},|x_0|\rangle+\frac{\alpha}{1-\alpha}\varpi(L^\alpha) 
\]
In this situation, we set $a_k :=  E[\|x_k\|^2_Q]$ to get from Proposition \ref{prop.tauberian} (ii) that
\[
\lim_{\alpha \to 1} (1 - \alpha) \sum_{k=0}^{\infty} \alpha^k E[\|x_k\|^2_Q]  
   =\lim_{\kappa\to\infty}\frac{1}{\kappa}\sum_{k=0}^{\kappa-1}  E[\|x_k\|^2_Q]=\varpi(L),
\]
since both limits exist, showing the result. 
\end{proof}

\subsection{Output Energy Measurements and $H_2$ Norms}

 For the discounted $H_2$-norm and counter-discounted energy measurements, consider particularly Corollaries \ref{coro:stability.discounted} and   \ref{coro:detectability.stability}, Lemma \ref{lemma:converegence.alfa.qqer} and Remark \ref{rem:convergence.properties}. 

\begin{corollary}(The discounted $H_2$-norm and the counter-discounted energy measurement) \label{coro:discounted.norm} Suppose that $\Theta$ is $(C, \mathcal{L}^\alpha)$-detectable and $L^\alpha \succeq 0$ is the solution of \eqref{lyapunovkind.1} with $Q=C^\trp C$. If $\alpha>1$ assume also that $r_\sigma(A)<\alpha^{-1}$. Then, system $\Theta$ is $\alpha$-stochastically stable and
\[
\lim_{\kappa \to \infty}\sum_{k=0}^{\kappa} E\Bigl[\alpha^{k}\left( \|y_k\|^2-\alpha\varpi(L^\alpha) \right) \Big]<\infty
\]
holds.  When $\alpha<1$, its mean energy norm is
\begin{equation}
    \mathcal{E}^\alpha_2 = \lim_{\kappa \to \infty}  E \Big[\sum_{k=0}^{\kappa} \alpha^k\| y(k) \|^2 \Big] = \frac{\alpha}{1-\alpha}\varpi(L^\alpha).
\end{equation}
When $\alpha\ge1$, there is $\beta>0$ such that $\bigl|E[ \|y_k\|^2]-\alpha\varpi(L^\alpha)\bigr]\bigr|\le \beta \alpha^{-k}$.
\end{corollary}

For the long-run $H_2$-norm, Theorem \ref{thm:analysis:longrunnorm2} completes the necessary elements.
\begin{corollary}($H_2$-norm: The long-run case) \label{coro:average.norm} Suppose that $\Theta$ is $(C, \mathcal{L})$-detectable and $L\succeq0$ is the solution of \eqref{convergent_lyapunov} with $Q=C^\trp C$. Then, system $\Theta$ is stochastically stable and its mean power norm is
\begin{equation}
    \mathcal{P}_2 = \lim_{\kappa \to \infty} \frac{1}{\kappa} E \Big[\sum_{k=0}^{\kappa} \| y(k) \|^2 \Big]     = \varpi(L)
\end{equation}
Moreover,  $E\bigl[ \|y_k\|^2 \bigr]\to \varpi(L)$ as $k\to\infty$ and 
\begin{equation}
    \lim_{\alpha \to 1, \kappa \to \infty}(1-\alpha)\sum_{k=0}^{\kappa} E\Bigl[\alpha^{k}\left( \|y_k\|^2-\alpha\varpi(L^\alpha) \right) \Big]
    \\
    =\mathcal{P}_2-\varpi(L)
\end{equation}
\end{corollary}

\section{Conclusion}
\label{sec:conclusion}

The paper presents the essentials of CSVIU modeling as a contribution to the uncertain systems literature. It develops the energy measurement notion as a $H_2$-norm in Abel's mean form, or as overtaking measurement comparison, in addition to the power $H_2$-norm analysis associated with the C\`{e}saro's mean. 

It establishes the connections between stochastic stability and finiteness of total energy measurements, power measurement, or bounded energy increase rate of a CSVIU system, depending on the chosen notion. These different criteria deal with distinct focuses on transient and long-run behavior. The concept of $\alpha$-detectability and appropriate notions of energy or power measurements allow us to link such measurements and stability. A solution to a modified Lyapunov equation is essential to stochastic stability but not sufficient for the counter-discounted problem explored here. 

The long-run problem and the discounted/counter-discounted problems are connected utilizing Tauberian theorems, made precise the notion of vanishing discount involving the corresponding Abel's and C\`{e}saro's mean. The task of relating the discounted/counter-discounted with $\alpha$ ``nearly arbitrary" and the long-run norm is thus completed, which amounts to set $\alpha=1$ in most of the statements.

The definition of the power norm and the counter-discounted measurements are essential to reinforce the system state convergence in the mean-square sense to a known value, which the discounted norm problem renders undetermined.  In addition, the counter-discounted measurement yields a geometric estimate on the convergence rate that becomes stricter as the $\alpha$  parameter increases.

\bibliography{Daniel}

\begin{thebibliography}{10}

\bibitem{bib:19}
{\sc A.~E. Bouhtouri, D.~Hinrichsen, and A.~J. Pritchard}, {\em
  $\mathcal{H}_\infty$-type control for discrete-time stochastic systems}, Int.
  J. Robust and Nonlinear Control, 9 (1999), pp.~923--948.

\bibitem{EFCosta}
{\sc E.~F. Costa and J.~B.~R. do~Val}, {\em On the observability and
  detectability of continuous-time {Markov} jump linear systems}, SIAM J.
  Control Optim., 41 (2002), pp.~1295--1314.

\bibitem{bib:2}
{\sc O.~L.~V. Costa, M.~D. Fragoso, and R.~P. Marques}, {\em Discrete-Time
  Markov Jump Linear Systems}, Springer-Verlag, London-UK, 2005.

\bibitem{bib:17}
{\sc Damm, P.~Benner, and J.~Hauth}, {\em Computing the stochastic
  $\mathcal{H}_\infty$-norm}, arXiv:1703.04440v1,  (2017).

\bibitem{Davis}
{\sc M.~H.~A. Davis and R.~B. Vinter}, {\em Stochastic Modelling and Control},
  Chapman and Hall, London, 1985.

\bibitem{CSVIU:control}
{\sc J.~B.~R. do~Val and D.~S. Campos}, {\em The ${H}_2$-optimal control
  problem of {CSVIU} systems: Discounted, counter-discounted and long-run
  solutions, part {II}: Optimal control}, SIAM J. Control Optim., submitted
  (2021).

\bibitem{bib:4}
{\sc J.~B.~R. do~Val and R.~F. Souto}, {\em Modeling and control of stochastic
  systems with poorly known dynamics}, IEEE Transactions on Automatic Control,
  62 (2017), pp.~4467--4482, \url{https://doi.org/10.1109/TAC.2017.2668359}.

\bibitem{DraganMorozanStoica}
{\sc V.~Dragan, T.~Morozan, and A.~Stoica}, {\em Robust Control Theory},
  Springer, 2006.

\bibitem{bib:11}
{\sc V.~Dragan, T.~Morozan, and A.~Stoica}, {\em Mathematical methods in robust
  control of linear stochastic systems}, Springer, New York, 2nd~ed., 2013.

\bibitem{Fernandes2020}
{\sc M.~R. Fernandes, J.~B.~R. do~Val, and R.~F. Souto}, {\em Robust estimation
  and filtering for poorly known models}, IEEE Control Systems Letters, 4
  (2020), pp.~474--479, \url{https://doi.org/10.1109/LCSYS.2019.2951611}.

\bibitem{FrancisKhar}
{\sc B.~A. Francis and P.~P. {Khargonekar~(eds.)}}, {\em Robust Control
  Theory}, The IMA Volumes in Mathematics and its Applications, Springer, 1995.

\bibitem{freiling}
{\sc G.~Freiling and A.~Hochhaus}, {\em Properties of the solutions of rational
  matrix difference equations.}, Computers \& Mathematics with Applications, 36
  (2003), pp.~1137--1154.

\bibitem{Hasanov}
{\sc V.~I. Hasanov}, {\em Perturbation theory for linearly perturbed algebraic
  {R}iccati equations}, Numerical Functional Analysis and Optimization, 35
  (2014), pp.~1532--1559, \url{https://doi.org/10.1080/01630563.2014.895765}.

\bibitem{bib:8}
{\sc O.~Hernández-Lerma and J.~B. Lasserre}, {\em Discrete-Time Markov Control
  Processes: Basic Optimality Criteria}, Springer-verlag, New York - USA, 1996.

\bibitem{bib:9}
{\sc O.~Hernández-Lerma and T.~Prieto-Rumeau}, {\em {The vanishing discount
  approach to average reward optimality: the strongly and the weakly continuous
  cases}}, Morfismos, 12 (2008), pp.~1--15.

\bibitem{bib:18}
{\sc D.~Hinrichsen and A.~J. Pritchard}, {\em Stochastic $\mathcal{H}_\infty$},
  SIAM J. Control Optim., 36 (1998), pp.~1504--1538.

\bibitem{MKac}
{\sc M.~Kac}, {\em Some mathematical models in science}, Science, New Series,
  166 (1969), pp.~695--699.

\bibitem{bib:1}
{\sc M.~Kac, G.-C. Rota, and J.~T. Schwartz}, {\em Discrete Thoughts: Essays on
  Mathematics, Science and Philosophy}, Springer, New York, USA, 1992.

\bibitem{Kailath}
{\sc T.~Kailath}, {\em Linear Systems}, Prentice-Hall, Inc., Englewood Cliffs -
  USA, 1980.

\bibitem{Khasminskii.2012}
{\sc R.~Khasminskii}, {\em Stochastic Stability of Differential Equations}, 2nd
  Edition, Springer, New York - USA, 2012.

\bibitem{Zhao-Yanetall}
{\sc Z.-Y. Li, Y.~Wang, B.~Zhou, and G.-R. Duan}, {\em Detectability and
  observability of discrete-time stochastic systems and their applications},
  Automatica, 45 (2009), pp.~1340--1346,
  \url{https://doi.org/10.1016/j.automatica.2009.01.014}.

\bibitem{Mao.1991}
{\sc X.~Mao}, {\em Stability of Stochastic Differential Equations with Respect
  to Semimartingales}, Longman Higher Education, University of Michigan - USA,
  1991.

\bibitem{MeynTweedie}
{\sc S.~P. Meyn and R.~L. Tweedie}, {\em Markov Chain and Stochastic
  Stability}, Springer, London, 1993.

\bibitem{bib:3}
{\sc F.~C. Pedrosa, J.~C. Nereu, and J.~B.~R. do~Val}, {\em When control and
  state variations increase uncertainty: Modeling and stochastic control in
  discrete time}, Automatica, 123 (2021), p.~109341,
  \url{https://doi.org/10.1016/j.automatica.2020.109341}.

\bibitem{bib:6}
{\sc R.~Sznadjer and J.~A. Filar}, {\em Some comments on a theorem of {Hardy}
  and {Littlewood}}, J. Optim. Theory Appl., 75 (1992), pp.~201--208.

\bibitem{Walters}
{\sc P.~Walters}, {\em An Introduction ot Ergodic Theory}, Springer, London,
  2005.

\end{thebibliography}

\end{document}